\documentclass[12pt, twoside]{amsart}

\newtheorem{theorem}{Theorem}

\newtheorem{prop}[theorem]{Proposition}

\usepackage{amssymb,latexsym, amsmath, amscd}

\def\po{\parindent 0pt}
\def\p{{\po \it \s Proof. }}

\def\s{\smallskip}
\def\m{\medskip}

\def\cat{\operatorname{cat}}

\def\eps{\varepsilon}

\begin{document}

\title{On sequences cat\,$(f^k)$}


\author{Yu. B. Rudyak}
\address{Yu.\ B.\ Rudyak, Department of Mathematics, University of Florida, 358 
Little Hall, PO Box 118105
Gainesville,  FL 32611-8105, USA} \email{rudyak@math.ufl.edu}

\begin{abstract} 
We prove that, for every decreasing sequence $\{a_k\}$ of natural numbers, 
there exists a map $f: X \to X$ with $\cat f^k=a_k$.
\end{abstract}

\subjclass {Primary 55M30}
\maketitle

Let $f: X \to Y$ be a map of finite $CW$-spaces. An {\it $f$-categorical set} is 
an open subset $A$ of $X$ such that $f|A: A \to Y$ is null-homotopic. An {\it 
$f$-categorical covering} is a covering $\{A_i\}$ of $X$ such that every set 
$A_i$ is $f$-categorical. The Lusternik--Schnirelmann  category $\cat f$ of $f$ 
is defined to be the minimal number $k$ such that there exists a $k$-elemented 
$f$-categorical covering, \cite{F, Fe, BG}. We also set $\cat X  
=\cat  1_X$.

Given a pointed space, we always denote its base point by $*$.

\begin{prop}[\cite{BG}]
\label{p1}
For every diagram 
$$
\CD
X @>f>> Y @>g>> Y.
\endCD
$$  
we have: $\cat (gf) \le \min\{\cat f, \cat g\}$ 
\end{prop}

\p Clearly, every $f$-categorical covering is $gf$-categorical one. 
Furthermore, 
if $\{A_i\}$ is a $g$-categorical covering then $\{f^{-1}(A_i)\}$ is 
$gf$-categorical covering.
\qed

\begin{prop}[cf.  \cite{F}]\label{p2}
For two maps $f: X \to Y$ and $g: A \to B$ we have
$$
\cat( f \vee g) = \max\{\cat f, \cat g\}
$$ 
where $f\vee g: X \vee A \to Y \vee B$ is the wedge of $f$ and $g$. 
\end{prop}

\p We assume that $\cat f =k$ and $\cat g =l$ with $k\ge l$. Let $\{U_1, \ldots, 
U_k\}$ be an $f$-categorical covering and let $\{V_1, 
\ldots, V_l\}$ be a $g$-categorical covering. The image of $U_i$ in 
$X\vee A$ will also be denoted by $U_i$, and similarly for $V_j$. We can and 
shall assume that $*\notin V_i$ for $i>1$ and $*\notin V_j$ for $j>1$. Then, 
clearly, $U_i\cup V_i, 1<i\le l$ is an $f\vee g$-categorical set. Furthermore, 
since
$$
(X\vee A)\setminus (U_1 \cup V_1) = (X\setminus U_1) \cup (A\setminus V_1),
$$
we conclude that $U_1\cup V_1$ is open in $X\vee A$, and therefore $U_1\cup 
V_1$ is $f\vee g $-categorical. Now,
$$
U_1\cup V_1, \ldots, U_l \cup V_l,\, U_{l+1}, \ldots, U_k
$$
is an $f\vee g$-categorical covering. 
\qed

\m
Consider now a map $f: X \to X$. Because of what we said above, the sequence 
$\{\cat (f^k)\}$ is a decreasing sequence of natural numbers. Here we prove 
that, conversely,  every such a sequence can be realized as the sequence of the 
form $\cat{f^k}$. In other words, the following theorem holds.

\m 
{\noindent \bf Theorem.} 
{\it Let $\{a_k\}$ be a sequence of natural numbers with $a_k \ge a_{k+1}$. 
Then  there exists a self map $f: X \to X$ of a $CW$-space such that 
$\cat(f^k)=a_k$ and $\cat X=\cat f=a_1$.}

\p Clearly, every sequence $\{a_k\}$ as above stabilizes, i.e. there exists $m$ 
such that $a_k=a_{k+1}$ for $k\ge m$. We set
$$
D(\{a_k\})=\min \{m \bigm | a_k=a_{k+1} \text{ for } k\ge m\}
$$
and prove the theorem by induction on $D(\{a_k\})$.

For $D=1$ the result is obvious. Namely, we take a space $X$ with $\cat X=a_1$ 
(for example, the real projective space of dimension $a_1-1$) and set $f=1_X$. 

Now, suppose that the theorem holds for $D=n$, i.e. that for every sequence  
$\{c_k\}$ with $D(\{c_k\})=n$ there exists a map $f: X \to X$ with 
$\cat(f^k)=c_k$. Consider a sequence $\{a_k\}$ with $D(\{a_k\})=n+1$ and set 
$b_k=a_{k+1}, k=1, 2, \ldots$. Then $D(\{b_k\})=n$, and so there exists a map 
$f: X \to X$ with $\cat X=b_1$ and $\cat (f^k)=b_k$. 

Let $Y$ be a pointed space with $\cat Y=a_1$, and let $\eps: Y \to \{*\} 
\subset Y$ be the constant map. Consider the map
$$
\CD
u: Y \vee Y @>1\vee \eps >> Y \vee Y @>\pi >> Y = \{*\} \vee Y \subset Y \vee Y
\endCD
$$
where $\pi=\pi_Y: Y \vee Y \to Y$ is the pinch map. Clearly, $\cat u\le \cat Y$. 
On the other hand, since the composition
$$
\CD
Y=Y\vee \{*\} \subset Y \vee Y @>u>> Y \vee Y @>\pi>> Y
\endCD
$$
is the identity map, we conclude that, by Proposition \ref{p1}, $\cat  Y =\cat 
1_Y \le \cat u$. Thus, $\cat u =\cat Y=a_1$.

Consider the map 
$$
\CD
v: X \vee X @>1 \vee f >> X \vee X @>\pi >> X = \{*\} \vee X \subset X \vee X
\endCD
$$
and set $g=u\vee v: Y \vee Y \vee X \vee X  \to Y \vee Y \vee X \vee X$. We 
claim that $\cat (g^k)=a_k$.

First, by Proposition \ref{p2}, for $k=1$ we have
$$
\cat g= \max\{\cat u, cat v\}=\max\{a_1, \cat v\}=a_1,
$$
because $\cat v \le \cat(X \vee X)=\cat X =b_1\le  a_1$.

Furthermore, we have $\cat g^k=\max \{\cat u^k, \cat v^k\}$. But $u^k$ is the 
constant map for $k>1$, and hence $\cat g^k=\cat v^k$. So, it suffices to prove 
that $\cat v^k=\cat f^{k-1}$ for $k>1$.

It is easy to see that $v^k$ has the form
$$
\CD
v^k: X \vee X @>f^{k-1}\vee f^k >> X \vee X @>\pi >> X = \{*\} \vee X \subset X 
\vee X.
\endCD
$$ 
So, by Propositions \ref{p1} and \ref{p2},
$$
\cat v^k \le \cat(f^{k-1} \vee f^k) =\cat f^{k-1}.
$$
On the other hand, $f^{k-1}$ can be decomposed as
$$
\CD
f^{k-1}: X =X \vee \{*\} \subset X \vee X @>v^k >> X \vee X @>\pi >> X,
\endCD
$$
and so $\cat f^{k-1} \le \cat v^k$. Thus, $\cat v^k =\cat f^{k-1}$.

\m Finally, $\cat (Y \vee Y \vee X \vee X)=\cat Y=a_1$ by Proposition \ref{p2}. 
This completes the proof. 
\qed

\m {\bf Acknowledgment.} I am grateful to Alexander Felshtyn for statement of 
the problem: he asked me about properties of sequences of the form $\cat f^k$.

\end{document}